\documentclass[11pt]{amsart}
\usepackage{a4wide}
\usepackage[latin1, utf8]{inputenc}
\usepackage{lmodern}

\usepackage[hang,flushmargin]{footmisc}

\usepackage{amsmath,amssymb,amsfonts}
\usepackage[foot]{amsaddr}
\usepackage{graphicx}
\usepackage{hyperref}
\usepackage{orcidlink}

\usepackage[square,numbers]{natbib}
\usepackage{doi}

\usepackage[font=footnotesize,labelfont=footnotesize]{caption}

\usepackage{xurl}
\usepackage{hyperref}

\newtheorem{theorem}{Theorem}[section]
\newtheorem{proposition}{Proposition}[section]

\usepackage{mathtools}
\mathtoolsset{showonlyrefs}
\author[M.~Faulhuber]{Markus Faulhuber \orcidlink{0000-0002-7576-5724}}
\address{Faculty of Mathematics, University of Vienna \newline Oskar-Morgenstern-Platz 1, 1090 Vienna, Austria}
\email{markus.faulhuber@univie.ac.at}

\begin{document}
\title[On a fundamental barrier of the Wirtinger criterion]{On a fundamental barrier of the Wirtinger criterion \\ for Gabor systems with odd functions}

\thanks{This research was funded in whole or in part by the Austrian Science Fund (FWF) [\href{https://doi.org/10.55776/PAT5102224}{10.55776/PAT5102224}].}



\begin{abstract}
    We show that the Wirtinger criterion cannot be used to investigate the frame set conjecture for the first Hermite function. More generally, for odd functions, it cannot determine regions of the frame set with density less than 2.
\end{abstract}

\maketitle

\vspace*{-2\baselineskip}
\section{Introduction}
We consider the Gabor system $\mathcal{G}_1(a,b) = \{ M_{b \ell} T_{a k} h_1 \mid (k,\ell) \in \mathbb{Z}^2\}$. This is a Gabor system over a rectangular lattice $a \mathbb{Z} \times b \mathbb{Z}$ with the window being the first Hermite function
\begin{equation}
    h_1 (t) = c_1 \, t \, e^{-\pi t^2}, \quad t \in \mathbb{R},
\end{equation}
where $c_1 \in \mathbb{C}\setminus\{0\}$ is an arbitrary constant, which we can ignore for our analysis. The product $ab = \delta$ is the co-volume of the rectangular lattice $a \mathbb{Z} \times b \mathbb{Z}$ in $\mathbb{R}^2$. The quantity $\delta^{-1} = (ab)^{-1}$ is the density of the lattice and gives the average number of lattice points per unit area. We formulate our results for the co-volume, avoiding the need to take reciprocals repeatedly. 

A result of Gröchenig and Lyubarskii \cite{GroechenigLyubarskii_Hermite_2007, GroechenigLyubarskii_SuperHermite_2009} implies that $\mathcal{G}_1$ is a frame for $L^2(\mathbb{R})$ once $ab<1/2$. In the other direction, we know by the results of Lyubarskii and Nes \cite{LyubarskiiNes_Rational_2013} that $\mathcal{G}_1$ cannot be a frame if $ab = n/(n+1)$, $n$ being any positive integer. The frame set conjecture for the first Hermite function (cf.\ \cite{Gro14}) states that $\mathcal{G}_1$ is a frame whenever $ab<1$ and $ab \neq n/(n+1)$. Yet, we lack further knowledge when $1/2<ab<1$, except for numerics carried out in \cite{LyubarskiiNes_Rational_2013} for $a=1$, supporting the conjecture, and the remark in \cite[\S~5]{LyubarskiiNes_Rational_2013} that $\mathcal{G}_1(a,b)$ is a frame for $ab = 3/5$, with the additional comment that the computations are cumbersome and become impassable for general rational values of the co-volume $ab$.

In \cite{FauShaZlo25}, the Janssen test was used to find new frames with Hermite functions of order 4 and higher, and it is remarked that the criterion does not shed any new light on the case $n=1$. In \cite{GhoSel25}, Ghosh and Selvan used the \textit{Wirtinger criterion}, which we present in \S~\ref{sec:Wirtinger}, to numerically study Gabor systems with Hermite functions (and other window functions). The method is based on the Wirtinger inequality, which is the reason why we call it the Wirtinger criterion. The Wirtinger inequality was used earlier
by Sun and Zhou \cite{SunZho03} to prove a stability result for irregular Gabor frames. The numerics in \cite{GhoSel25} did not provide new insights into the frame set of $h_1$. This note aims to show that their criterion can indeed not be used to enlarge the frame set of $h_1$. To the best of our knowledge, no impossibility results have been established for sufficiency criteria in this context. Our contribution shows a fundamental obstacle intrinsic to the method, and that there is no further benefit in applying the criterion to the frame set conjecture for $h_1$.

\section{Preliminaries}
We follow the common notation in time-frequency analysis, as in the book of Gröchenig~\cite{Gro01}. We denote a time-frequency shift by $\pi(z) = M_\omega T_x$, $z = (x,\omega) \in \mathbb{R}^2$. Here, $T_x f(t) = f(t-x)$ is a translation (or time-shift) and $M_\omega f(t) = e^{2 \pi i \omega t} f(t)$ is a modulation (frequency-shift). A Gabor system $\mathcal{G}(g,\Gamma) = \{\pi(\gamma) g \mid \gamma \in \Gamma \subset \mathbb{R}^2 \text{ discrete}\}$ is a frame for $L^2(\mathbb{R})$ if there exist positive constants $0 < A \leq B < \infty$ such that
\begin{equation}\label{eq:frame}
    A \|f\|_{L^2}^2 \leq \sum_{\gamma \in \Gamma} |\langle f, \pi(\gamma) g \rangle|^2 \leq B \|f\|_{L^2}^2, \quad \forall f \in L^2(\mathbb{R).}
\end{equation}
In this case any $f \in L^2(\mathbb{R})$ can be written as a series of the form
\begin{equation}
    f = \sum_{\gamma \in \Gamma} c_\gamma \, \pi(\gamma) g, \quad (c_\gamma)_{\gamma \in \Gamma} \in \ell^2(\Gamma).
\end{equation}
We refer to \cite{Christensen_2016} for more details on frames. A function is suitable if it is in the Wiener space of continuous functions $W_0(\mathbb{R})$ with its Fourier transform also in $W_0(\mathbb{R})$. The Fourier transform of a suitable function $g$ will be denoted by $\widehat{g}$ and is given by
\begin{equation}
    \widehat{g}(\omega) = \int_\mathbb{R} g(t) e^{-2 \pi i \omega \cdot t} \, dt.
\end{equation}
Moreover, we use the unitary dilation operator, which is given by
\begin{equation}\label{eq:dilation_operator}
    \mathcal{D}_a f(t) = a^{-1/2} f(t/a), \quad a > 0.
\end{equation}
It is straightforward to check (see, e.g., \cite{Fau25, Fol89, Gos15, Gro01}) the following intertwining relation:
\begin{equation}\label{eq:dilation}
    \mathcal{D}_a \, \pi(z) \, \mathcal{D}_a^{-1} = \pi(D_a z), \quad
    D_a = \begin{pmatrix}
        a & 0 \\
        0 & 1/a
    \end{pmatrix}.
\end{equation}
Based on this result, one can easily prove the following result (see, e.g., \cite{Fau25, Gos15, Gro01}).
\begin{proposition}\label{pro:unitary}
    The Gabor system $\mathcal{G}(g,a\mathbb{Z} \times b \mathbb{Z)}$ is unitarily equivalent to the Gabor system $\mathcal{G}(\mathcal{D}_b g, \delta\mathbb{Z} \times \mathbb{Z})$, $\delta = ab$. In particular, one system is a frame if and only if the other is.
\end{proposition}

\section{The Wirtinger criterion}\label{sec:Wirtinger}
We recall a sufficiency criterion for Gabor systems, which we call the Wirtinger criterion. It is based on the following Wirtinger inequality with (at least) one boundary point fixed to~0.
\begin{equation}
    \int_0^1 |f(t)|^2 \, dt \leq \frac{4}{\pi^2}\int_0^1|f'(t)|^2 \, dt, \quad f(0) \cdot f(1)=0, \ f \in \mathcal{C}^1([0,1]).
\end{equation}
The inequality can be used to derive a sufficient criterion for a shift-invariant system to form a frame. This criterion can be transferred to Gabor systems. We refer to \cite{GhoSel25, GroRomSto18, SunZho03} for details.

The criterion below has been formulated in \cite[Thm.~2.4]{GhoSel25}, where also the technical assumptions (suitability, differentiability, decay) on $g$ are clarified in detail.
\begin{theorem}\label{thm:Wirtinger}
    The Gabor system $\mathcal{G}(g, \delta \mathbb{Z} \times \mathbb{Z}$) is a frame if
    \begin{equation}
        \delta < \delta_g(\omega) = \frac{1}{2} \sqrt{\frac{\sum \limits_{k \in \mathbb{Z}} |\widehat{g}(k+\omega)|^2}{\sum \limits_{k \in \mathbb{Z}} (k+\omega)^2 |\widehat{g}(k+\omega)|^2}}, \quad \forall \omega \in [0,1].
    \end{equation}
\end{theorem}

\section{The Wirtinger criterion and the Gaussian window}
For illustrative purposes, we first show that the criterion almost gives the optimal bound for the Gaussian Gabor system $\mathcal{G}(\varphi, \delta \mathbb{Z} \times \mathbb{Z})$, $\varphi(t) = c_0 \, e^{- \pi t^2}$, $c_0 \in \mathbb{C} \setminus \{0\}$. Note that $\widehat{\varphi} = \varphi$.

As shown in \cite[Lem.~4.2]{GhoSel25}, the function
\begin{equation}
    \omega \mapsto \frac{\sum \limits_{k \in \mathbb{Z}} |\varphi(k+\omega)|^2}{\sum \limits_{k \in \mathbb{Z}} (k+\omega)^2 |\varphi(k+\omega)|^2}, \quad \omega \in [0,1],
\end{equation}
attains its minimum if $\omega = 1/2$. We now bound the numerator and denominator. We have
\begin{equation}
    \sum_{k \in \mathbb{Z}} |\varphi(k+1/2)|^2 > \sum_{k\in\{-1,0\}} |\varphi(k+1/2)|^2 = 2 e^{-\pi/2}.
\end{equation}
In the other direction, it is easy to see that (note that $x \mapsto x e^{-\pi x}$ is decreasing for $x>1/\pi$)
\begin{align}
    \sum_{k \in \mathbb{Z}} (k+1/2)^2 |\varphi(k+1/2)|^2 & = 2 \sum_{k \geq 0} (k+1/2)^2 |\varphi(k+1/2)|^2\\
    & < \frac{e^{-\pi/2}}{2} + 2 \sum_{k \geq 1} (k+1/2) e^{-2 \pi (k+1/2)}
\end{align}
The exponential tails can be easily bounded from above. We have
\begin{align}
    \sum_{k\geq 1} e^{-2 \pi (k+1/2)} & = \frac{e^{-\pi}}{e^{2\pi}-1} \approx 0.00008085 < 10^{-4},\\
    \sum_{k\geq 1} k \, e^{-2 \pi (k+1/2)} & = \frac{e^{\pi}}{(e^{2\pi}-1)^2} \approx 0.000081 < 10^{-4}.
\end{align}
Thus, we have that $2 \sum_{k \geq 1} (k+1/2) e^{-2 \pi (k+1/2)} < 3 \cdot 10^{-4}$ and
\begin{equation}
    \frac{\sum \limits_{k \in \mathbb{Z}} |\varphi(k+\omega)|^2}{\sum \limits_{k \in \mathbb{Z}} (k+\omega)^2 |\varphi(k+\omega)|^2} > 4 \, \frac{e^{-\pi/2}}{e^{-\pi/2} + 6 \cdot 10^{-4}} \approx 4 \cdot 0.997122.
\end{equation}
Note that we could have been more careful in our analysis, but already with these bounds, the Wirtinger criterion yields that $\mathcal{G}(\varphi, \delta \mathbb{Z} \times \mathbb{Z})$ is a frame once
\begin{equation}
    \delta \leq 0.9985 < \frac{1}{2} \sqrt{4 \cdot 0.997112} < \delta_\varphi(\omega), \quad \forall \omega \in [0,1].
\end{equation}
By employing Proposition~\ref{pro:unitary} and similar estimates as above, we can analyze Gabor systems of the form $\mathcal{G}( \mathcal{D}_b \varphi, \delta \mathbb{Z} \times \mathbb{Z})$ which are unitarily equivalent to $\mathcal{G}(\varphi, a \mathbb{Z} \times b \mathbb{Z})$, $\delta = ab$. As illustrated in \cite[Fig.~2(A)]{GhoSel25}, one can get close to the truth ($ab < 1$) for the frame set of the Gaussian function $\varphi$, by using Theorem~\ref{thm:Wirtinger}. As we shall see next, for $h_1$ and, in fact, any odd function, the Wirtinger criterion cannot overcome the fundamental barrier $ab = 1/2$.

\section{The Wirtinger criterion and \texorpdfstring{$h_1$}{h1}}
We consider the dilated first Hermite function
\begin{equation}
    \phi_b(t) = \mathcal{D}_b h_1(t) = \frac{c_1}{\sqrt{b}} \frac{t}{b} e^{-\frac{\pi}{b^2} t^2}, \quad t \in \mathbb{R}, \ b >0.
\end{equation}
Note that it satisfies the functional equation $\widehat{\phi_b} = -i \, \phi_{1/b}$. Next, we consider the function
\begin{equation}
    \delta_{\phi_b}(\omega)
    = \frac{1}{2} \sqrt{\frac{\sum \limits_{k \in \mathbb{Z}} |\phi_{1/b}(k+\omega)|^2}{\sum \limits_{k \in \mathbb{Z}} (k+\omega)^2 |\phi_{1/b}(k+\omega)|^2}} 
    = \frac{1}{2} \sqrt{\frac{\sum \limits_{k \in \mathbb{Z}} (k+\omega)^2 e^{-2\pi b^2(k+\omega)^2}}{\sum \limits_{k \in \mathbb{Z}} (k+\omega)^4 e^{-2\pi b^2(k+\omega)^2}}}, \quad b > 0.
\end{equation}
Now, if $\delta < \min_{\omega} \delta_{\phi_b}$, then the Gabor system $\mathcal{G}(\phi_b, \delta \mathbb{Z} \times \mathbb{Z})$ is a frame. This is equivalent, by Proposition~\ref{pro:unitary}, to the fact that $\mathcal{G}_1(a,b)$ is a frame, $ab = \delta < \min_{\omega} \delta_{\phi_b}(\omega)$.

In order to gain new insights into the frame set of $h_1$, we would need to find parameters $b \in \mathbb{R}_+$ such that $\min_{\omega}\delta_{\phi_b}(\omega) \geq 1/2$. We prove that this is not possible, implying that the Wirtinger criterion cannot be used to further study the frame set of $h_1$.

As can be easily seen, we have
\begin{equation}
    \min_{\omega} \delta_{\phi_b}(\omega) \leq \delta_b(0) = \frac{1}{2} \sqrt{\frac{\sum \limits_{k \in \mathbb{Z}} k^2 e^{-2 \pi b k^2}}{\sum \limits_{k \in \mathbb{Z}} k^4 e^{-2 \pi b k^2}}} < \frac{1}{2}, \quad \forall b >0.
\end{equation}
This follows, as for each individual $k \in \mathbb{Z}\setminus\{0\}$, the corresponding term in the numerator is smaller than (or equal to) the corresponding term in the denominator. Generally,
\begin{equation}
    |\widehat{g}(k)|^2 \leq k^2 | \widehat{g}(k)|^2, \quad k \in \mathbb{Z} \setminus\{0\},
\end{equation}
with equality if and only if $k^2 = 1$. The only term that could make $\delta_g(0)$ large must come from $|\widehat{g}(0)|^2$, which equals zero for the first Hermite function.

In order to determine whether a Gabor system with the first Hermite function is a frame, the Wirtinger criterion will only become conclusive once the condition $ab < 1/2$ is satisfied. Thus, no selection of parameters $(a,b)$ within the Wirtinger criterion can resolve any part of the still‑open region $1/2 < a b < 1$ for the first Hermite function.

\section{Odd window functions}\label{sec:odd}
We now show that the Wirtinger criterion cannot determine points in the frame set of a suitable odd function if $1/2 < ab$. Recall that parameters with $ab=1/2$ are always excluded from the frame set of odd functions \cite{LyubarskiiNes_Rational_2013}. An extension to general lattices is given in \cite{Faulhuber_Note_2018}.

The argument is again the same as above. For an odd function $g(t) = -g(-t)$, we naturally have $g(0)=0$. The Fourier transform preserves parity, so $\widehat{g}(0)=0$. It follows that
\begin{equation}
    \sum_{k \in \mathbb{Z}} |\widehat{g}(k)|^2 < \sum_{k \in \mathbb{Z}} k^2 |\widehat{g}(k)|^2.
\end{equation}
This implies that, for any suitable odd $g$, we have
\begin{equation}
    \min_{\omega} \delta_g(\omega) \leq \delta_g(0) = \frac{1}{2} \sqrt{\frac{\sum \limits_{k \in \mathbb{Z}} |\widehat{g}(k)|^2}{\sum \limits_{k \in \mathbb{Z}} k^2 |\widehat{g}(k)|^2}} < 1/2.
\end{equation}
Thus, the Wirtinger criterion is inconclusive for any Gabor system of the form $\mathcal{G}(g, \delta \mathbb{Z} \times \mathbb{Z})$, $1/2 \leq \delta$. By Proposition~\ref{pro:unitary} this Gabor system is unitarily equivalent to $\mathcal{G}(\mathcal{D}_b^{-1} g, a \mathbb{Z} \times b \mathbb{Z})$, $ab = \delta$. As the dilation operator preserves a function's parity, the limitation follows.

\section{General lattices}
The above barrier cannot be bypassed by general (non-rectangular) lattices either. This follows from the following facts. Any lattice $\Lambda$ in $\mathbb{R}^2$ is symplectic up to scaling, i.e.,
\begin{equation}
    \Lambda = \delta S \mathbb{Z}^2, \quad S \in SL(2,\mathbb{R}) = Sp(1), \ \delta > 0.
\end{equation}
We refer the reader to the textbook by Folland for details \cite{Fol89} on the symplectic group $Sp(d) \subset SL(2d,\mathbb{R})$. By an Iwasawa decomposition, we can write
\begin{equation}
    \Lambda = \delta \ R_r V_q D_a \mathbb{Z}^2,
    \quad \text{ where } \quad
    R_r = \begin{pmatrix}
        \cos(r) & \sin(r)\\
        -\sin(r) & \cos (r)
    \end{pmatrix},
    \quad
    V_q = \begin{pmatrix}
        1 & 0\\
        q & 1
    \end{pmatrix},
    \quad r, \ q \in \mathbb{R}
\end{equation}
and $D_a$ was defined in \eqref{eq:dilation}. The two-fold cover of the symplectic group is the metaplectic group \cite{Gos11}. The corresponding metaplectic operators are the dilation operator defined in~\eqref{eq:dilation_operator},
\begin{equation}
    \mathcal{F}_r f(s) = \int_{\mathbb{R}} f(t) \ k_r(s,t) \, dt,
    \quad
    \mathcal{V}_q f(t) = e^{\pi i q t^2} \ f(t),
\end{equation}
with $k_r(s,t) = \sqrt{1-i \cot(r)} \ \exp(\pi i (\cot(r) s^2 - 2 \csc(r) s t + \cot(r) t^2)$. The operators $\mathcal{F}_r$ and $\mathcal{V}_q$ are the fractional Fourier transform and the chirp. We note that metaplectic operators preserve the parity of functions \cite[Chap.~4.4]{Fol89} (cf.\ \cite{Fau19_RiM}). Proposition~\ref{pro:unitary} indeed holds more generally. Assume $\mathcal{U}$ is a unitary operator on $L^2(\mathbb{R})$ and $U$ a symplectic matrix in $Sp(1)$ which together satisfy
\begin{equation}\label{eq:unitary}
    \mathcal{U} \, \pi(z) \, \mathcal{U}^{-1} = c_U(z) \, \pi(U z), \quad |c_U(z)|=1.
\end{equation}
Then, see \cite{Fau25_SampTA} (cf.\ \cite{Gos15}, \cite[Chap.~9.4]{Gro01}), the Gabor systems $\mathcal{G}(g,\Lambda)$ and $\mathcal{G}(\mathcal{U}g , U \Lambda)$ are unitarily equivalent. We provide the necessary computations from \cite{Fau25_SampTA}. Denote the frame operator of the Gabor system $\mathcal{G}(g,\Lambda)$ by
\begin{equation}
    S_{\mathcal{G}(g,\Lambda)} = \sum_{\lambda \in \Lambda} \langle f, \pi(\lambda) g \rangle \pi(\lambda) g.
\end{equation}
Then $A$ and $B$ from \eqref{eq:frame} are its spectral bounds. We have
\begin{align}
    \mathcal{U} \, S_{\mathcal{G}(g,\Lambda)} \, \mathcal{U}^{-1} f
    = & \ \sum_{\lambda \in \Lambda} \langle \mathcal{U}^{-1} f, \pi(\lambda) g \rangle \, \mathcal{U} \pi(\lambda) g
    = \sum_{\lambda \in \Lambda} \langle f, \mathcal{U} \pi(\lambda) \mathcal{U}^{-1} \mathcal{U} g \rangle \, \mathcal{U} \pi(\lambda) \mathcal{U}^{-1} \mathcal{U} g\\
    = & \ \sum_{\lambda \in \Lambda} \langle f, \pi(U \lambda) \mathcal{U} g \rangle \, \pi(U \lambda) \mathcal{U} g
    = S_{\mathcal{G}(\mathcal{U}g, U \Lambda)} f
\end{align}
Note that in the above calculations the phase factor $c_U$, appearing in \eqref{eq:unitary} also appears as complex conjugate $\overline{c_U}$ and, so, cancels. Equation \eqref{eq:unitary} holds for any metaplectic operator and its projection onto the symplectic group, in particular for the pairs $(\mathcal{D}_a, D_a)$, $(\mathcal{V}_q, V_q)$, and $(\mathcal{F}_r, R_r)$ (see \cite{Fau25_SampTA} and the references therein). As the parity of functions is kept under the metaplectic deformations, we can reduce our analysis to the rectangular (indeed square) lattice case in \S~\ref{sec:odd}. Thus, the Wirtinger criterion cannot detect frames with odd functions and $1/2 \leq \text{vol}(\mathbb{R}^2/\Lambda)$.


\begin{thebibliography}{99}
\renewcommand*{\doi}[1]{\href{https://doi.org/#1}{#1}}

\bibitem{Christensen_2016}
O.~Christensen.
\newblock \emph{{An Introduction to Frames and Riesz Bases}}.
\newblock {Applied and Numerical Harmonic Analysis}. Birkhäuser, 2. edition, 2016.
\newblock \doi{10.1007/978-3-319-25613-9}

\bibitem{Faulhuber_Note_2018}
M.~Faulhuber.
\newblock {A short note on the frame set of odd functions}.
\newblock \emph{Bulletin of the Australian Mathematical Society}, 98(3):481--493, 2018.
\newblock \doi{10.1017/S0004972718000746}.

\bibitem{Fau19_RiM}
M.~Faulhuber.
\newblock {On the Parity under Metaplectic Operators and an Extension of a Result of Lyubarskii and Nes}.
\newblock \emph{Results in Mathematics}, 75, 2020.
\newblock \doi{10.1007/s00025-019-1134-4}

\bibitem{Fau25}
M.~Faulhuber.
\newblock {Gabor systems with Hermite functions of order $n$ and oversampling greater than $n+1$ which are not frames}.
\newblock \emph{{Sampling Theory, Signal Processing, and Data Analysis}}, 23:22, 2025.
\newblock \doi{10.1007/s43670-025-00112-5}.

\bibitem{Fau25_SampTA}
M.~Faulhuber.
\newblock {Operator-isomorphism pairs and Zak transform methods for the study of Gabor systems}.
\newblock In \emph{{Proceedings of the International Conference on Sampling Theory and Applications}}, 2025.
\newblock \doi{10.1109/SampTA64769.2025.11133554}.

\bibitem{FauShaZlo25}
M.~Faulhuber, I.~Shafkulovska, and I.~Zlotnikov.
\newblock {On the frame property of Hermite functions and explorations of their frame sets}.
\newblock \emph{{Journal of Fourier Analysis and Applications}}, 31(21), 2025.
\newblock \doi{10.1007/s00041-025-10153-y}.

\bibitem{Fol89}
G.~B. Folland.
\newblock \emph{{Harmonic Analysis in Phase Space}}.
\newblock Number 122 in {Annals of Mathematics Studies}. Princeton University Press, 1989.
\newblock \doi{10.1515/9781400882427}.

\bibitem{GhoSel25}
R.~Ghosh and A.~A. Selvan.
\newblock {On Gabor frames generated by B-splines, totally positive functions, and Hermite functions}.
\newblock \emph{Applied Numerical Mathematics}, 207:1--23, 2025.
\newblock \doi{10.1016/j.apnum.2024.08.021}.

\bibitem{Gos11}
M.~A.\ de~Gosson.
\newblock \emph{{Symplectic Methods in Harmonic Analysis and in Mathematical Physics}}, volume~7 of \emph{{Pseudo-Differential Operators. Theory and Applications}}.
\newblock Birkhäuser/Springer Basel AG, Basel, 2011.
\newblock \doi{10.1007/978-3-7643-9992-4}.

\bibitem{Gos15}
M.~A.\ de~Gosson.
\newblock {Hamiltonian deformations of Gabor frames: First steps}.
\newblock \emph{Applied and Computational Harmonic Analysis}, 38(2):196--221, 2015.
\newblock \doi{10.1016/j.acha.2014.03.010}.

\bibitem{Gro01}
K.~Gröchenig.
\newblock \emph{{Foundations of Time-Frequency Analysis}}.
\newblock {Applied and Numerical Harmonic Analysis}. Birkhäuser, Boston, MA, 2001.
\newblock \doi{10.1007/978-1-4612-0003-1}.

\bibitem{Gro14}
K.~Gröchenig.
\newblock {The Mystery of Gabor Frames}.
\newblock \emph{Journal of Fourier Analysis and Applications}, 20(4):865--895, 2014.
\newblock \doi{10.1007/s00041-014-9336-3}.

\bibitem{GroechenigLyubarskii_Hermite_2007}
K.~Gröchenig and Y.~Lyubarskii.
\newblock {Gabor frames with Hermite functions}.
\newblock \emph{Comptes Rendus Mathematique}, 344(3):157--162, 2007.
\newblock \doi{10.1016/j.crma.2006.12.013}.

\bibitem{GroechenigLyubarskii_SuperHermite_2009}
K.~Gröchenig and Y.~Lyubarskii.
\newblock {Gabor (super)frames with Hermite functions}.
\newblock \emph{Mathematische Annalen}, 345(2):267--286, 2009.
\newblock \doi{10.1007/s00208-009-0350-8}.

\bibitem{GroRomSto18}
K.~Gröchenig, J.~L. Romero, and J.~Stöckler.
\newblock {Sampling theorems for shift-invariant spaces, Gabor frames, and totally positive functions}.
\newblock \emph{Inventiones mathematicae}, 211(3):1119--1148, 2018.
\newblock \doi{10.1007/s00222-017-0760-2}.

\bibitem{LyubarskiiNes_Rational_2013}
Y.~Lyubarskii and P.~G. Nes.
\newblock {Gabor frames with rational density}.
\newblock \emph{Applied and Computational Harmonic Analysis}, 34(3):488--494, 2013.
\newblock \doi{10.1016/j.acha.2012.09.001}.

\bibitem[Sun and Zhou(2003)]{SunZho03}
W.~Sun and X.~Zhou.
\newblock {Irregular Gabor frames and their stability}.
\newblock \emph{Proceedings of the American Mathematical Society}, 131:2883--2893, 2003.
\newblock \doi{10.1090/S0002-9939-02-06931-9}.

\end{thebibliography}
\end{document}